\documentclass[12pt]{article}

\usepackage{times}
\usepackage[a4paper,text={128mm,185mm}, centering]{geometry}
\usepackage{changepage}
\usepackage{titlesec}

\titleformat{\section}[hang]%
{\bfseries\large}{\thesection.}{1ex}{}%

\titleformat{\subsection}[hang]%
{\bfseries}{\thesubsection}{1ex}{}%

\usepackage{amsmath}
\usepackage{fancyhdr}
\usepackage{amsthm, amssymb, bm}
\usepackage[all]{xy}

\title{\vskip 5pt  \bf  THE TOPOLOGY OF CRITICAL PROCESSES,  III   $\hspace{-100pt}$
\newline \rm (COMPUTING HOMOTOPY)}
\author{\itshape\bfseries {Marco GRANDIS}}
\date{}

\def \LL {\skp \setlength{\leftskip}{8pt}}  %left indent, for all the following
\def \LB {\setlength{\leftskip}{0pt}}   %end of left indent
\let \lan \langle
\let \ran \rangle
\def \sst {\scriptstyle}

\def \c {\colon}

\def \q {\qquad}
\def \qq {\qquad \qquad}

\def \bu {{\scriptscriptstyle\bullet}}
\def \skp {\medskip}
\def \ndt {\noindent}
\def \Ndt {\medskip  \noindent}
\def\shp {^{\sharp}}
\def \tilde {\raise.17ex\hbox{$\scriptstyle\mathtt{\sim}$}}   %for web address
\def \ul {\underline}

\def \uw {\!\, \raisebox{0.3ex} {$\uparrow$}}

\def \ard {\ar@{-->}}
\def \arp {\ar@{.>}}
\def \are {\ar@{->>}}
\def \aru {\ar@{=}}
\def \arv {\ar@{}}   %void arrow, to put a label
\def \arl {\ar@{-}}    %line
\def \arld {\ar@{--}}    %dashed line
\def \arlp {\ar@{..}}    %dotted line

\def \ti {\! \times \!}

\def \me {{\, {\scriptstyle{\wedge}}\, }}
\def \ci {{\raise.3ex\hbox{$  \scriptscriptstyle\circ $}}}

\def \setm {{\raise.4ex\hbox{$ \, \scriptscriptstyle{\setminus} \; $}}}
\def \sub {\subset}
\def\le{\leqslant}
\def\ge{\geqslant}
\def \eqq {\! \sim \!_2 \,}

\def \and {\mbox{ and }}

\def \id {{\rm id\,}}

\def \al {\alpha}

\def \ph {\varphi}

\def \one {{\bf 1}}
\def \two {{\bf 2}}
\def \three {{\bf 3}}
\def \N {{\bf N}}   %order category
\def \Z {{\bf Z}}   %order category
\def \R {{\bf R}}   %order category
\def\bc {{\bf c}}

\def \Cat {\mathsf{Cat}}

\def \dTop {{\rm d}\mathsf{Top}}
\def \cTop {{\rm c}\mathsf{Top}}

\def \cbfTop {{\rm c_{bf}}\mathsf{Top}}

\def \Fl {\mathsf{F}{\rm l}\,}
\def \uPi {\uw\Pi}
\def \upi {\uw\pi}

\def \bbN {{\mathbb{N}}}  % {{...}}  for exponentiation
\def \bbZ {\mathbb{Z}}
\def \bbI {{\mathbb{I}}}   
\def \bbR {{\mathbb{R}}}
\def \bbS {\mathbb{S}}

\def \uI {{\uw\mathbb{I}}}
\def \uR {\uw\mathbb{R}}
\def \uS {\uw\mathbb{S}}

\def \rc {{\rm c}}
\def \cI {{\rm c}\mathbb{I}}
\def \cJ {{\rm c}\mathbb{J}}
\def \cR {{\rm c}\mathbb{R}}
\def \cS {{\rm c}\mathbb{S}}
\def \cT {{\rm c}\mathbb{T}}

\begin{document}
\maketitle

\vskip 25pt
% Abstract followed by Keywords and MSC:
\begin{adjustwidth}{0.5cm}{0.5cm}
{\small

{\bf Abstract.} Directed Algebraic Topology studies spaces equipped with a form of direction, to include 
models of non-reversible processes. In the present extension we also want to cover {\em critical 
processes}, indecomposable and unstoppable.

	The previous parts of this series introduced {\em controlled spaces} and their fundamental category. 
Here we study how to compute the latter. The homotopy structure of these spaces will be examined in 
Part IV.\\


}
\end{adjustwidth}

% Here starts the text

%%%
%%%
\section*{Introduction}\label{Intro}

%%%
\subsection{Directed and controlled spaces}\label{0.1}
Directed Algebraic Topology is an extension of Algebraic Topology, dealing with `spaces' where the paths 
need not be reversible; the general aim is including the representation of {\em irreversible processes}. A 
typical setting for this study, the category $\dTop$ of directed spaces, or d-spaces, was introduced and 
studied in \cite{G1}--\cite{G3}; it is often employed in the theory of concurrency, cf.\ \cite{FGHMR}.

	The present series is devoted to a further extension, where the paths can also be non-decomposable 
in order to include {\em critical processes}, indivisible and unstoppable -- either reversible or not. 
For instance: quantum effects, the onset of a nerve impulse, the combustion of fuel in a piston, the switch 
of a thermostat, the change of state in a memory cell, the action of a siphon, moving in a no-stop road, etc.

	To this effect the category of d-spaces was extended in Part I  \cite{G4} to the category $\cTop$ of 
{\em controlled spaces}, or {\em c-spaces}: an object is a topological space equipped with a set $X\shp$ 
of continuous mappings $a\c [0, 1] \to X$, called {\em controlled paths}, or {\em c-paths}, which are closed 
under concatenation and global reparametrisation (by surjective increasing endomaps of the interval) 
and include all the constant paths at the endpoints of c-paths.

	A {\em map of c-spaces}, or {\em c-map}, is a continuos mapping which preserves the selected paths. 
Their category $\cTop$ contains the category $\dTop$ of d-spaces as a full subcategory, reflective and 
coreflective: a c-space is a d-space if and only if it is {\em flexible}, which means that each point is flexible 
(its trivial loop is controlled) and every controlled path is flexible (all its restrictions are controlled).

	Every c-space $X$ has two associated d-spaces, the generated d-space $\hat{X}$ and the flexible 
part $\Fl X$, by the reflector and coreflector of the embedding $\dTop \to \cTop$ (Section 1.2 of Part I).

%%%
\subsection{The fundamental category}\label{0.2}
	Part II \cite{G5} defines and studies the fundamental category of controlled spaces, as a functor
    \begin{equation}
\uPi_1 \c \cTop \to \Cat,
    \label{0.2.1} \end{equation}
that extends the fundamental category of d-spaces \cite{G1, G3} and the fundamental groupoid of 
topological spaces. 

	There are two natural transformations (see Section 5.2 of Part II)
    \begin{equation}
\uPi_1(\Fl X) \, \longrightarrow \, \uPi_1(X)  \, \longrightarrow \, \uPi_1(\hat{X})
    \label{0.2.2} \end{equation}
induced by the embeddings $\Fl X \to X \to \hat{X}$ (the counit of the coreflector and the unit of the reflector 
of d-spaces). 

	These functors need not be faithful, as we shall see in \ref{1.3}, but Theorem 5.3(b) of Part II says that 
$\uPi_1(X) \to \uPi_1(\hat{X})$  is a full embedding when the c-space $X$ is {\em preflexible}, that is
all the c-paths of $\hat{X}$ between flexible points of $X$ are already controlled in the latter.

	The present Part III is an immediate continuation of Part II, devoted to computing the fundamental 
category of c-spaces. The definitions and results of Part II are taken for granted and only referred to.
	
	Part IV will study the homotopy structure of c-spaces, their homotopy equivalences and their links 
with cubical sets. In particular, we shall analyse the formal theory of homotopy in $\cTop$, following 
the classification of directed settings in \cite{G3}.

%%%
\subsection{Outline}\label{0.3}
In Section 1 we calculate the fundamental category of the c-spaces introduced so far, and others, applying 
Theorems 5.3 (on preflexible c-spaces) and 5.8 (on covering maps of c-spaces) of Part II, and developing 
peculiar techniques adequate to the present framework. The relationship between the fundamental category  
of c-spaces and d-spaces is discussed in \ref{1.6}, where we show that the theorem of Seifert-van Kampen 
fails for c-spaces.

	In the same line, Section 2 briefly considers how the analysis of obstructions, a typical problem in 
concurrency, can be dealt with replacing the d-spaces used in \cite{G3}, Chapter 3 (and elsewhere) with 
rigid c-spaces. This leads to a far simpler analysis, but a less rich one.

	Finally, in Section 3, we prove that the fundamental category of a border flexible c-space can be simply 
defined by general deformations of controlled paths, instead of using their flexible deformations -- as in the 
general case.

\Ndt {\em Acknowledgments}. The author is indepted to the Referee for many helpful suggestions.

%%%
\subsection{Notation and conventions}\label{0.4}
A continuous mapping between topological spaces is called a {\em map}. $\bbR $ denotes the euclidean 
line as a topological space, and $\bbI$ the standard euclidean interval $[0, 1]$. The identity path $\id \bbI$ 
is written as $\ul{i}$. The open and semiopen intervals of the real line are denoted by square brackets, 
like $]0, 1[$, $[0, 1[$ etc.

	A {\em preorder} relation is assumed to be reflexive and transitive; an {\em order} is also 
anti-symmetric. A mapping which preserves (resp.\ reverses) preorders is said to be {\em increasing} 
(resp.\ {\em decreasing}), always used in the weak sense. 

	As usual, a preordered set $ X $ is identified with the small category whose objects are the elements 
of $ X$, with one arrow $ x \to x' $ when $x$ precedes $x'$ and none otherwise.
 
	The binary variable $\al$ takes values $ 0, 1, $ which are generally written as $ -, + $ in superscripts 
and subscripts. The symbol $\, \sub \,$ denotes weak inclusion.

\vskip 5pt
	The previous papers \cite{G4, G5} of this series are cited as Part I and Part II, respectively; the reference 
I.2 or II.3.4, for instance, points to Section 2 of Part I or Subsection 3.4 of Part II.

%%%
%%%
\section{Calculating the fundamental category}\label{s1}

	This section studies how to compute the fundamental category of c-spaces. Using Theorem II.5.3(b) 
on preflexible c-spaces, many of these results can be deduced from the fundamental category of the 
generated d-spaces, already computed in \cite{G3}; but a direct calculation can often be simple and 
more significant.

	The new aspects which appear here, with respect to the theory of d-spaces, are highlighted in \ref{1.6}.

\begin{small}
\vskip 5pt

	The symbols $\two, \three, \N, \Z, \R$ denote ordered sets, and the associated categories; the ordered 
sets  2,  3  and  $D|\Z|$ are discrete. $\bbN$ is the one-object category associated to the additive monoid of 
the natural numbers.

\end{small}

%%%
\subsection{Elementary calculations}\label{1.1}
We begin by examining the basic c-spaces, showing that many of them are 1-simple, in the sense of 
II.5.1: their fundamental category is a preorder; of course, the controlled circles $\cS^1$ and $\rc_n\bbS^1$ 
are not. (Some of these results are already in II.5.9.)

\Ndt (a) The fundamental categories of $ \cI$, $ \cJ$, $ \cR $ are the following ordered sets:
    \begin{equation}
\uPi_1(\cI)  =  \two,    \q   \uPi_1(\cJ)  =  \three,    \q   \uPi_1(\cR)  =  \Z.
    \label{1.1.1} \end{equation}

	As to $\cI$, the identity $ \ul{i}\c \cI \to \cI $ is 2-equivalent to any other c-path $ \rho\c 0 \to 1$, by 
Lemma II.4.6(c): in fact, $ \rho $ is a global reparametrisation, and therefore $ \rho = \ul{i}\rho \eqq \ul{i}$, 
so that there is precisely one arrow $ [\ul{i}] $ from 0 to 1, in the fundamental category. At each flexible point, 
0 or 1, there is only one loop $ \cI \to \cI$, the trivial one.

\begin{small}
\vskip 3pt

	As to $\cJ$ and $\cR$, two c-paths $ a, b\c x \to y$ in any of them are always 2-equivalent, since they are 
in the one-jump c-structure of $[x, y]$, isomorphic to $\cI$. 

\end{small}
\vskip 3pt

	For these preflexible spaces the components of the natural transformations  
$\uPi_1(\Fl X) \to \uPi_1(X) \to \uPi_1(\hat{X})$ of \eqref{0.2.2} become inclusions of ordered sets:
    \begin{equation}
2 \to \two \to [0, 1],    \q   3 \to \three \to [0, 2],    \q   D|\Z| \to \Z \to \R.
    \label{1.1.2} \end{equation}

\Ndt (b) The argument used above for $ \uPi_1(\cI) $ also applies to the delayed intervals $\rc_-\bbI$ and 
$\rc_+\bbI$, in II.1.3(e)
    \begin{equation}
\uPi_1(\rc_-\bbI)  =  \uPi_1(\rc_+\bbI)  =  \two,
    \label{1.1.3} \end{equation}
whose c-structure is also generated by a single map $ \bbI \to \bbI$. These c-spaces are not preflexible, but 
their fundamental category is still full in $ \uPi_1(\uI)$.

\Ndt (c) The fundamental category of the directed circle $ \uS^1$, as described in \cite{G3}, 3.2.7(d), is the 
subcategory of the groupoid $ \Pi_1 \bbS^1 $ formed of the classes of anticlockwise paths (in $ \bbR^2$). 
Each monoid $ \upi_1(\uS^1, x) $ is isomorphic to the additive monoid $\bbN$ of natural numbers.

	Applying Theorem II.5.3(b), the fundamental category of the one-stop circle $ \cS^1 $ amounts to the 
fundamental monoid at the unique flexible point $x_0$ (the point 1 of the complex plane)
    \begin{equation}
\uPi_1(\cS^1)(x_0, x_0)  =  \upi_1(\uS^1, x_0)  =  \bbN.
    \label{1.1.4} \end{equation}

	Without using $\uPi_1(\uS^1)$ this is also proved by Theorem II.5.8(b) applied to the exponential map 
$ \cR \to \cS^1$.

	Therefore two c-loops $ a, b $ in $ \cS^1 $ are 2-equivalent if and only if they have the same length 
$ 2k\pi$ (in radians), if and only if they both turn $ k $ times ($k \ge 0$) around the circle, anticlockwise.

\Ndt (d) More generally, the fundamental category of the preflexible $n$-stop circle $ \rc_n\bbS^1 $ (see 
II.1.4(d) is the full subcategory of the fundamental category of $(\rc_n\bbS^1)\hat{\,} = \uS^1 = \uR/\bbZ$ 
on $ n $ flexible points, the vertices $ [i/n] $ (for $ i = 0, ..., n-1$) of an inscribed $n$-gon.

	$\uPi_1(\rc_n\bbS^1)$ is thus the category $\bc_n $ freely generated by $ n $ arrows disposed as 
follows on the edges of an $n$-gon
%
%  FIG. 1.1.5: the n-gon categories
    \begin{equation} 
\xy <1pt, 0pt>:
%  bullets, labels, dummy pt above, below
(20,0) *{\bu}; (70,0) *{\bu}; (110,0) *{\bu}; 
(200,0) *{\bu}; (170,17) *{\bu}; (170,-17) *{\bu}; 
(35,-12) *{\bc_1}; (125,-12) *{\bc_2}; (215,-12) *{\bc_3}; 
(0,25) *{}; (0,-25) *{};
% arrows
(20, 5); (-20, 0) **\crv{(18,20)&(-20,20)}, (20, -5); (-20, 0) **\crv{(18,-20)&(-20,-20)}, 
(110, 4); (70, 4) **\crv{(105,20)&(75,20)}, (110, -4); (70, -4) **\crv{(105,-20)&(75,-20)}, 
\POS(19.3, -7) \ar+(1,3), \POS(71.5,7) \ar-(1,3), \POS(108.5, -7) \ar+(1,3), 
\POS(195,3) \ar+(-20,12), \POS(170,12) \ar-(0,23), \POS(175,-15) \ar+(20,12), 
\endxy
    \label{1.1.5} \end{equation}

	Again, this result can also be obtained using the covering map of c-spaces 
$p_n\c \rc_n\bbR \to \rc_n\bbS$.

\Ndt (e) For the preflexible c-space $ X $ on the euclidean interval $ [0, 3] $ described in I.2.3(e) we have a 
mixed situation; essentially, the paths in $ [1, 2] $ behave as in $ \cI$, while those in $ [0, 1] $ or $ [2, 3] $ 
behave as in $ \uI$.

%%%
\subsection{Higher dimensional c-spaces}\label{1.2}
(a) Applying Theorem II.5.6 on cartesian products, we get the following fundamental categories
    \begin{equation} \begin{array}{c}
\uPi_1(\cI^n)  =  \two^n,   \q   \uPi_1(\cJ^n)  =  \three^n,
\\[7pt]
\uPi_1(\cI \ti \cJ)  =   \two \ti \three,	
\\[5pt]
\uPi_1(\cR^n)  =  \Z ^n,   \qq   \uPi_1(\cT^n)  =  \bbN^n,
    \label{1.2.1} \end{array} \end{equation}
which are (partially) ordered sets, except the last. The controlled $n$-torus $\cT^n$ was defined in I.2.6(d) 
as the cartesian power $(\cS^1)^n$, or equivalently as the orbit c-space $ (\cR^n)/\bbZ^n$; its fundamental 
category amounts to the monoid $\bbN^n$ at the only flexible point.

\Ndt (b) The fundamental category of all the higher c-spheres $ \cS^n$, for $ n \ge 2$, is trivial: the discrete 
category $\one$.

	In fact, there is one flexible point, $*$. Every c-path of $ \cS^n $ is a general concatenation of a finite 
family of c-loops of the form $pa$, where $ a\c \cI \to \cI^n $ is a c-path of the controlled $n$-cube, and it is 
sufficient to prove that each of them is 2-equivalent to the trivial loop (at $*$).

	If the path $ a $ lies in a face of the cube, $ pa $ is already the trivial loop. Otherwise, it is a path 
$(0, \, ..., 0) \to (1, \, ..., 1)$, 
and it is 2-equivalent to the concatenation $ b = b_1 * b_2 $ of two c-paths living in some faces, and collapsed 
to the trivial loop in the quotient c-space. 	For instance one can take $ b_1(t) = (t, 0, ..., 0) $ (on an edge) and 
$ b_2(t) = (1, t, ..., t) $ (in the face $ t_1 =1$).

%%%
\subsection{Other calculations}\label{1.3}
The following computations of the fundamental category give a better understanding of the natural 
transformations $\, \uPi_1(\Fl X) \to \uPi_1(X) \to \uPi_1(\hat{X})$ of \eqref{0.2.2}. Moreover, they are based 
on topological arguments which will also be useful in other cases.

\Ndt (a) The reversible c-interval $ \cI^\sim $ of II.1.3(d) has a c-structure generated by the identity path 
$ \ul{i} $ and the reversion $ r\c \bbI \to \bbI$; the flexible points are 0 and 1.

	Each c-path $ x \to y $ (between flexible points) has an integral length, which is even if $ x = y $ and odd if 
$ x \neq y$. We prove below, in Theorem \ref{1.7}, that this length is constant up to 2-equivalence, and 
determines the class of a path in $ \uPi_1(\cI^\sim)(x, y)$.

	In other words, we shall prove that the obvious c-map $ p\c \rc_2\bbS^1 \to \cI^\sim $
%
%  FIG. 1.3.1
    \begin{equation} 
\xy <1pt, 0pt>:
%  bullets, labels, dummy pt above, below
(-20,0) *{\bu}; (20,0) *{\bu}; (-20,-32) *{\bu}; (20,-32) *{\bu}; 
(110,0) *{p(x, y) = (x + 1)/2}; (0,25) *{}; (0,-32) *{};
%
% circle, line
(0,0) = "A"  *\cir<20pt>{}="ca"  
\POS(-20,-32)\arl+(40,0)
% arrows
\POS(0,20)\are-(4,0), \POS(0,-20)\are+(4,0), \POS(-8,-32)\are-(4,0), \POS(8,-32)\are+(4,0), 
\endxy
    \label{1.3.1} \end{equation}
induces an isomorphism $ p_*\c \uPi_1(\rc_2\bbS^1) \to \uPi_1(\cI^\sim)$ defined on the category 
$ \bc_2 $ described in \eqref{1.1.5}. Let us note that $ p $ is not a covering map: the flexible points of the 
basis are not evenly covered; loosely speaking, the selection of c-paths in the domain and codomain 
`mends' this failure.

	Thus the category $\uPi_1(\cI^\sim)$ is freely generated by two arrows, the classes $ [\ul{i}]\c 0 \to 1 $ 
and $ [r]\c 1 \to 0$; at each vertex it has a fundamental monoid isomorphic to the additive monoid $ \bbN$.

	The generated d-space $(\rc_S\bbI)\hat{\,} = \bbI^\sim$ is the reversible d-interval of I.2.4(c), whose 
fundamental category is the indiscrete groupoid on two objects (with one arrow between any pair of objects). 

	In this case the functor $\uPi_1(X) \to \uPi_1(\hat{X})$ is not faithful; moreover $\ul{i}$ 
and $ r $ are reversible c-paths of $X$ whose classes in $\uPi_1(X)$ are not invertible.

\Ndt (b) The fundamental category $ \uPi_1(\rc_S\bbI) $ of the growing-siphon interval (in I.3.3(a)) is 
generated by the following arrows (where $ r $ is the reversion path $ r(t) = 1 - t$)
    \begin{equation}
(x, x')\c x \to x',   \q   [r]\c 1 \to 0   \qq   (0 \le x < x' \le 1),
    \label{1.3.2} \end{equation}
under the relation $ (x, x')(x', x'') = (x, x'')$, for $ 0 \le x < x' < x'' \le 1$.

	The identity path $ \ul{i} $ is flexible and reversible in $\rc_S\bbI$, but is not flexibly reversible: the 
reversed path $ r $ is not flexible, and the associated arrow $ [\ul{i}] = (0, 1)\c 0 \to 1 $ is not invertible. But 
it becomes invertible in the fundamental category of $ (\rc_S\bbI)\hat{\,} = \bbI^\sim$: also here the 
functor $\uPi_1(X) \to \uPi_1(\hat{X})$ is not faithful.

%%%
\subsection{On-off controller}\label{1.4}
We now examine the c-space $ X $ built in I.3.1(a) to model an on-off controller (e.g.\ a thermostat) that 
oversees a variable $T$ (e.g.\ the temperature), counteracting its rising
%
% FIGURE 1.4.1, from I.3.1.1  On-off controller
    
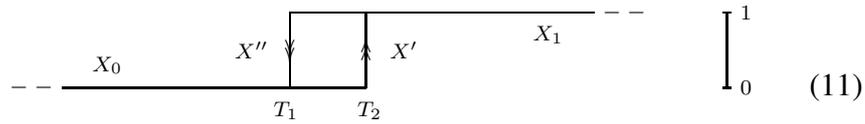
\begin{equation}
\xy <.5mm, 0mm>:
% labels + dummy pts left/above, right/below
(12,6) *{\sst{X_0}}; (50,10) *{\sst{X''}}; (90,10) *{\sst{X'}};  (128,14) *{\sst{X_1}}; 
(59,-6) *{\sst{T_1}}; (81,-6) *{\sst{T_2}}; 
(180,0) *{\sst{0}}; (180,20) *{\sst{1}};  
(-30,25) *{}; (180,-10) *{};
% branches X_0, X_1
\POS(0,0) \arl+(80,0),  \POS(60,20) \arl+(80,0),
\POS(-13,0) \arld+(10,0),  \POS(143,20) \arld+(10,0), 
% branches X', X'' and double arrows
\POS(60,0) \arl+(0,20),  \POS(80,0) \arl+(0,20),
\POS(60,12) \are+(0,-5),  \POS(80,8) \are+(0,5),
%vertical axis
\POS(175,0) \arl+(0,20), \POS(174,0) \arl+(2,0), \POS(174,20) \arl+(2,0), 
\endxy
    \label{1.4.1} \end{equation}

	On the left branch $ X_0 $ the system is in state 0: the cooling device is off; if the temperature grows to 
$ T_2 $ the device jumps to state 1; then, if the temperature cools to $ T_1$, it goes back to state 0.

	The support $ |X| $ of our model is a one-dimensional subspace of $ \bbR^2$. The c-structure of 
$ X $ is generated by the c-structures of:

\LL

\ndt - $ X_0$, $ X_1$, natural intervals where $ T $ can vary,

\Ndt - $ X'$, $ X''$, one-jump c-intervals, where $ T $ is constant and the state of the system varies.

\LB
\skp

	The flexible part $ X_0 + X_1 $ of the c-space $ X $ is the sum of two natural intervals; its fundamental 
groupoid $ \Pi_1 (\Fl X) $ is the sum of the indiscrete groupoids on the same sets, categorically equivalent to 
the discrete groupoid $ 2 = \{0, 1\}$.

	The fundamental category $ \uPi_1(X) $ is equivalent to its skeleton, the full subcategory on two points 
$ x_0 \in X_0 $ and $ x_1 \in X_1$; the latter is isomorphic to the category $\bc_2 $ (see \eqref{1.1.5}).

%%%
\subsection{Transport networks and labelled graphs}\label{1.5}
Transport networks are usually modelled in graph theory, in an effective way as far as they do not interact 
with continuous variation. They can also be modelled by c-spaces, which allows us to combine them with 
planar or three-dimensional regions, as we have discussed in I.3.4.

	The fundamental category can be readily used to study such models. Controlled spaces can thus 
unify aspects of continuous and discrete mathematics, interacting with hybrid control systems and 
others sectors of Control Theory [Br, He].

%%%
\subsection{Comments}\label{1.6}
(a) The main method of calculation of the fundamental category for complex spaces, the theorem of 
Seifert-van Kampen, holds true in $\dTop$, in the fundamental-category version of \cite{G3}, 3.2.6, but 
fails here.

	For instance, we have seen that the category $ \uPi_1(\cI) = \two $ has one arrow $ 0 \to 1$. Now we can 
cover $ \cI $ with the open subspaces $ U = [0, 1[ $ and $V = \; ]0, 1]$, which only inherit the trivial loops at 0 and 
1, respectively. Their fundamental category has only these trivial arrows, and the pushout over $\uPi_1(U \cap V)$ 
(the empty category) gives the discrete category $ 2$.
 
\Ndt (b) Nevertheless, we have seen that the fundamental category $\uPi_1(X)$ of a rigid or `partially rigid' 
 c-space can be rather easy to compute without this tool -- or using it on $\uPi_1(\hat{X})$ when the original 
 c-space is preflexible.
 
\Ndt (c) In many cases $\uPi_1(X)$ is very small and easy to analyse, while $\uPi_1(\hat{X})$ gives a finer 
description, at the price of a complex analysis where the equivalence of categories is totally ineffective. This 
will show even more clearly in the next section.

%%%
\subsection{Theorem}\label{1.7}
{\em
The projection $ p\c \rc_2\bbS^1 \to \cI^\sim $ defined in \eqref{1.3.1} induces an isomorphism of categories 
$p_*\c \uPi_1(\rc_2\bbS^1) \to \uPi_1(\cI^\sim)$.
}
\begin{proof}
(a) The functor $ p_* $ is bijective on the objects, the flexible points. It is also full, because 
$ p\c \rc_2\bbS^1 \to \cI^\sim $ obviously satisfies the path-lifting property II.5.7(i) within c-paths: every 
c-path $ b\c y \to y' $ in $ \cI^\sim $ has a lifting $ a\c x \to x' $ in $ \cI^\sim$, determined by the starting point 
$ x \in F_y $ (unique in the present case).

	The length of $ b $ is an integer, equal to the length of $ a $ measured in half-circles.

	To prove that $ p_* $ is faithful we shall show that two c-paths $ b, b'\c y_0 \to y_1 $ in $ \cI^\sim $ which 
are 2-equivalent {\em have the same length}, so that any pair of their liftings in $ \rc_2\bbS^1 $ starting at the 
same point are also 2-equivalent; in other words one can lift along $ p $ the 2-equivalence relation -- if not the 
actual 2-paths. 

	For the sake of simplicity we suppose that $ y_0 = 0$, the case $ y_0 = 1 $ being similar. We 
use the path spaces $P(\bbI) = \bbI^{\bbI} \,$ and $P(\bbI^2)$ with the compact-open topology, determined by 
the metric $ d(c, c') = \max_t \, d(c(t), c'(t)) $ (and the euclidean metric on $ \bbI $ and $ \bbI^2$).

\Ndt (b) Let $ P_n $ be the subspace of $ P(\bbI) $ formed of the c-paths $ \cI \to \cI^\sim $ starting at $ 0$, 
of length $ n$; let $ P $ be their (disjoint) union. We prove now that each $ P_n $ is open in $ P$. 
(This amounts to saying that the length function $ P \to \bbN $ is continuous, which is not obvious as it fails 
on the whole path space $ \bbI^\bbI$.)

	It will be sufficient to show that any two c-paths $ b, b'\c 0 \to y $ with $ d(b, b') < 1/2 $ have the same length. 
If $ b $ has length $ n$, it determines a partition of the interval $ \bbI $ in $ n $ subintervals
    \begin{equation} \begin{array}{c}
0 = t_0  <  t_1  <  ...  <  t_n = 1,
\\[5pt]
b(t_0)  =  0,   \q   b(t_1)  =  1, \, ...   \q\;\;   b(t_n)  =  (1 - (-1)^n)/2,
    \label{1.7.1} \end{array} \end{equation}
and is {\em properly} increasing on $ [0, t_1]$, properly decreasing on $ [t_1, t_2]$, and so on (by `properly' 
we mean that it is not constant). There are $ n - 1 $ `inversions of monotony' (each of them occurring on a 
maximal closed subinterval where $ b $ is constant at 1 or 0, alternatively).

	The other path $ b'$, of length $ n'$, has $ b'(t_0) = 0 $ and $ b'(t_1) > 1/2$; because of the form of c-paths 
in $ \cI^\sim$, it must be properly increasing on some (at least one) subinterval of $ [0, t_1]$. It also has 
$ b'(t_2) < 1/2$, and must be properly decreasing on some subinterval of $ [t_1, t_2]$; and so on. Finally, it has 
at least as many inversions of monotony as $ b$, and $ n' \ge n$. By symmetry, $ n = n'$.

\Ndt (c) Let $ K\c \cI \times \uI \to \cI^\sim $ be a hybrid 2-path between the c-paths $ b, b'\c 0 \to y$. Proving that 
they have the same length will achieve the argument.

	The family of c-paths
    \begin{equation} 
u_t\c \cI \to \cI^\sim \ti \uI,   \q   u_t(s)  =  (s, t)   \qq   (t \in \bbI),
    \label{1.7.2} \end{equation}
gives an isometry $ u\c \bbI \to P(\bbI^2)$
$$
d(u_t, u_{t'})  =  {\rm max}_s \, d((s, t), (s, t'))  =  |t - t'|.
$$

	Composing $ u $ with the map $ K_*\c P(\bbI^2) \to P(\bbI) $ we get a continuous mapping
    \begin{equation} 
Ku\c \bbI \to P(\bbI),    \q    t  \mapsto  K_t = K(-, t)\c \bbI \to \bbI,
    \label{1.7.3} \end{equation}
whose values $ K_t $ are the intermediate c-paths of $ K $ (see II.4.4(a)). They belong to $ P$. Since 
$ Ku $ is defined on a connected space, all of them belong to the same subset $ P_n$, including $b$ and 
$ b'$.
\end{proof}
%

%%%
%%%
\section{Analysing obstructions}\label{s2}

	The analysis of obstructions inside a cubical directed space is a typical problem in concurrency, dealt 
with in \cite{FGHMR} and many papers (see Part I). It is also studied in \cite{G3}, Chapter 3, working 
with d-spaces. The corresponding problem in rigid c-spaces seems to be far simpler, although it can give 
a less fine analysis, as shown in \ref{2.3}.

%%%
\subsection{An elementary case}\label{2.1}
We begin with the `square annulus' $ X \sub \cI^2 $ represented below, namely the compact subspace of 
the standard c-square which is the complement of the open square $ ]1/3, 2/3[^2 $ (marked with a cross); 
the latter should be viewed as a single obstruction in an unstoppable process
%
% Fig 2.1.1
\begin{equation} \begin{array}{c}
\xy <.4mm, 0mm>:
%HOR positions: [0,40], [100,140], 
%
% labels, dummy pts below, above
(55,-10) *{X}; (170,-10) *{\uPi_1(X)};
(100,-20) *{0}; (100,20) *{x}; (140, -21) *{y}; (140, 20) *{1}; 
(0,-20) *{}; (0,25) *{};
% x in hole and non-comm. sq.
(20,0) *{\sst{\times}}; (120,0) *{\times}; 
% bullets of X
@i @={(0,-20), (40,-20), (40,20), (0,20)} @@{*{\bu}};
% space X, ext. boundary and hole
@i @={(0,-20), (40,-20), (40,20), (0, 20)},
s0="prev"  @@{;"prev";**@{-}="prev"}
@i @={(13,-7), (27,-7), (27,7), (13, 7)},
s0="prev"  @@{;"prev";**@{-}="prev"},
% arrows of the fund. cat.
\POS(108,-20)\ar+(24,0), \POS(108,20)\ar+(24,0), 
\POS(100,-12)\ar+(0,24), \POS(140,-12)\ar+(0,24), 
\endxy
\label{2.1.1} \end{array} \end{equation}

	Typically, in the analysis of concurrent processes, the obstruction represents a resource (e.g.\ a memory 
storage, an application, a printer) that two (or more) concurrent automata cannot engage at the same time. A 
path below {\em or} above the obstruction corresponds to priority of one of them. Modifying the picture, one 
can represent in a similar way an island in a stream or a one-dimensional obstacle in space-time, as in the 
Introduction to \cite{G3}.

\skp

	The fundamental category $\uPi_1(X)$ is represented in the right  diagramabove: it is 
generated by four arrows forming a non-commutative square, and has two arrows $0 \to 1$ (not 
drawn in the figure).

	Applying Theorem II.5.3(b) one can deduce this fact from the fundamental category of the generated 
d-space $ \hat{X} \sub \uI^2$, determined in \cite{G3}, 3.1.1. But a direct proof is rather simple. 

\begin{small}
\vskip 5pt

In fact, every c-path $ a\c 0 \to 1 $ in $ X $ meets the vertical strip 
$$ S = \; ]1/3, 2/3[ \, \times \, \bbI$$ 
in one connected component of 
$ S \cap X$, either below or above the obstruction. Suppose that $ a $ meets the lower component 
$ U = \; ]1/3, 2/3[ \times [0, 1/3[ $ (open in $ X$). The preimage $a^{-1}(U)$ is an open subinterval of $]0, 1[$ 
(by continuity and monotony), and we can suppose it is precisely $]1/3, 2/3[$, up to invertible reparametrisation 
and 2-equivalence. For a second path $ a' $ of the same kind and similarly reparametrised, we can suppose 
that $ a(t) \le a'(t) $ for $ t \in \bbI $ (replacing $ a $ with $ a \me a'$). 

	Now the affine interpolation $H$ from $a$ to $a'$ is a hybrid 2-path in $\cI^2$ and takes the interval 
$ ]1/3, 2/3[ $ to the rectangle $ U $ (by monotony), proving that $ a \eqq a' $ in $ X$. Similarly, two paths 
above the obstruction are 2-equivalent in $ X$. Finally, a c-path below the obstruction and another above 
are not even 2-equivalent in the underlying topological space.

\end{small}

%%%
\subsection{Two obstructions}\label{2.2}
We examine now two subspaces $Y, Z \sub \cI^2$ which arise from two obstructions, either appearing 
together (with respect to the generated path order, see I.1.8(c)) or one after the other. 

In both cases a direct 
computation is easy, if more complex than in the previous case; alternatively, one can deduce our results 
from the fundamental category of the generated d-spaces, described in \cite{G3}, 3.9.2 and 3.9.4(b).

\Ndt (a) {\em Simultaneous obstructions}. The first case can be modelled with the subspace $Y$ of $\cI^2$ 
represented below
%
% Fig 2.2.1, from DAT 39.2.1 modified
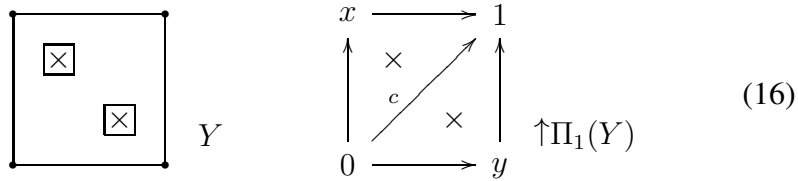
\begin{equation} \begin{array}{c}
\xy <.4mm, 0mm>:
%HOR POSITIONS: [10,60], [120,170],
%
% labels, dummy pts below, above
(75,20) *{Y}; (198,20) *{\uPi_1(Y)}; 
(120,10) *{0}; (120,60) *{x}; (170, 9) *{y}; (170, 60) *{1}; (135, 32) *{\sst{c}}; 
(0,3) *{};(0,67) *{};
% x in holes and non-comm. squares
@i @={(25,45), (45,25), (135,45), (155,25)} @@{*{\times}};
% bullets
@i @={(10,10), (60,10), (60,60), (10, 60)} @@{*{\bu}};
% space Y, ext boundary and two holes
@i @={(10,10), (60,10), (60,60), (10, 60)},
s0="prev"  @@{;"prev";**@{-}="prev"};
@i @={(20,40), (30,40), (30,50), (20, 50)},
s0="prev"  @@{;"prev";**@{-}="prev"};
@i @={(40,20), (50,20), (50,30), (40, 30)},
s0="prev"  @@{;"prev";**@{-}="prev"};
% arrows of the fund. cat.
\POS(128,10)\ar+(34,0), \POS(128,18)\ar+(34,34), \POS(128,60)\ar+(34,0), 
\POS(120,18)\ar+(0,34),  \POS(170,18)\ar+(0,34),
\endxy
\label{2.2.1} \end{array} \end{equation}

	The fundamental category $ \uPi_1(Y) $ has again four vertices; from 0 to 1 there are three arrows: 
$[a]$ (through $x$), $[b]$ (through $y$) and $[c]$.

\Ndt (b) {\em Consecutive obstructions}. The second case is modelled by $Z \sub \cI^2$
%
% Fig 2.2.2
    
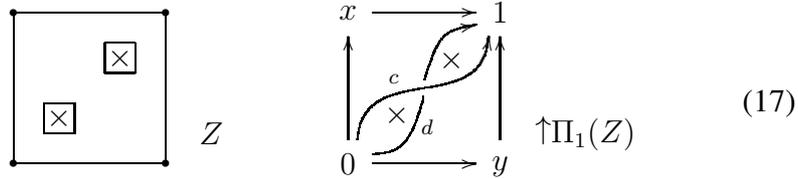
\begin{equation}    %%% one cannot use  \begin{array}  with  \crv !!!
\xy <.4mm, 0mm>:
%HOR POSITIONS: [10,60], [120,170],
%
% labels, dummy pts below, above
(75,-10) *{Z}; (198,-10) *{\uPi_1(Z)}; 
(120,-20) *{0}; (120,30) *{x}; (170, -21) *{y}; (170, 30) *{1}; (135, 8) *{\sst{c}}; (146, -8) *{\sst{d}}; 
(0,-25) *{};(0,40) *{};
% x in holes and non-comm. squares
@i @={(45,15), (25,-5), (154,14), (136,-4)} @@{*{\times}};
% bullets
@i @={(10,-20), (60,-20), (60,30), (10, 30)} @@{*{\bu}};
% space Y, ext boundary and two holes
@i @={(10,-20), (60,-20), (60,30), (10, 30)},
s0="prev"  @@{;"prev";**@{-}="prev"};
@i @={(40,10), (50,10), (50,20), (40, 20)},
s0="prev"  @@{;"prev";**@{-}="prev"};
@i @={(20,-10), (30,-10), (30,0), (20, 0)},
s0="prev"  @@{;"prev";**@{-}="prev"};
%% curved arrows
(123,-12); (166,22) **\crv{(125,15)&(165,-5)}, 
%(128,-17); (162,26) **\crv{(155,-15)&(135,25)}, 
(128,-17); (144.8,2) **\crv{(135,-16)&(140,-17)}, 
(145.2,8); (162,26) **\crv{(150,24)&(152,24)},
% arrows of the fund. cat. 
\POS(128,-20)\ar+(34,0), \POS(128,30)\ar+(34,0), 
\POS(120,-12)\ar+(0,34),  \POS(170,-12)\ar+(0,34),
\POS(166,19)\ar+(.5,3), \POS(158.7,25.2)\ar+(3,.5),
\endxy
\label{2.2.2} \end{equation}

	In $ \uPi_1(Z) $ there are now four arrows from 0 to 1: $[a]$ (through $x$), $[b]$ (through $y$) and $[c], [d]$.
	
\Ndt (c) {\em Comments}. The fundamental category distinguishes these situations, which topology cannot 
separate: the underlying topological spaces $ |Y| $ and $ |Z| $ are homeomorphic.

%%%
\subsection{Obstructions in d-spaces}\label{2.3}
The d-spaces $ \hat{X}, \hat{Y}, \hat{Z} $ generated by the previous c-spaces have the same topological 
support and the structure induced by the ordered square $ \uI^2$. 

	Their fundamental category, much more complex than in the previous cases, was studied in 
\cite{G3}, 3.1.1, 3.9.2, 3.9.4(b).

	In each case the fundamental category, whose objects are the infinite points of the support, is skeletal 
and {\em cannot} be reduced up to equivalence of categories. As analysed in \cite{G3}, Section 3.9, it is 
essentially represented by a `minimal injective model', future and past equivalent to the given category.
Here we get the finite, full subcategories represented below (on 4, 8, 6 objects, respectively), determining 
the `branching points' of the process
%
%
%  Fig 2.3.1,   from  DAT  31.1.2,  39.2.1,  39.4.3
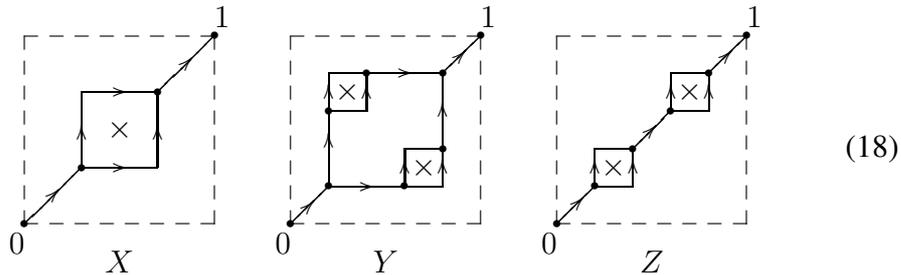
\begin{equation} \begin{array}{c}
\xy <.5 mm, 0mm>:
%POSITIONS: X [10,60],  Y [80,130], Z [150,200],
%
% labels, dummy pts below, above/right
(35,0) *{X}; (105,0) *{Y}; (175,0) *{Z}; 
(10,-8) *{}; (170,50) *{};
% x in holes
@i @={(35,35), (95,45), (115,25), (165,25), (185,45)} @@{*{\times}};
% bullets
@i @={(10,10), (25,25), (45,45), (60,60), 
(80,10), (90,20), (90,40), (100,50),(120,50), (130,60), (120,30), (110,20),
(150,10), (160,20), (170,30), (180,40), (190,50), (200,60)} @@{*{\bu}};
% labels
(8,5) *{0}; (62,65) *{1};    (78,5) *{0}; (132,65) *{1};    (148,5) *{0}; (202,65) *{1}; 
%  ext boundary
@i @={(10,10), (60,10), (60,60), (10, 60)},
s0="prev"  @@{;"prev";**@{--}="prev"},
@i @={(80,10), (130,10), (130,60), (80, 60)},
s0="prev"  @@{;"prev";**@{--}="prev"},
@i @={(150,10), (200,10), (200,60), (150, 60)},
s0="prev"  @@{;"prev";**@{--}="prev"}
%  paths from 0 to 1 (covering the bound. of the holes)
@i @={(10,10), (25,25), (25,45), (45, 45), (60, 60), (45, 45), (45, 25), (25, 25)},
s0="prev"  @@{;"prev";**@{-}="prev"},
@i @={(80,10), (90,20), (90,50), (120, 50), (130, 60),
(120, 50), (120, 20), (90, 20)},
s0="prev"  @@{;"prev";**@{-}="prev"}
@i @={(90,20), (90,40), (100,40), (100,50), (120, 50),
(120, 30), (110, 30), (110, 20)},
s0="prev"  @@{;"prev";**@{-}="prev"},
@i @={(150,10), (160,20), (160,30), (170,30), (180,40), (180, 50),
(190, 50), (200, 60), (190, 50), (190, 40), (180, 40), (170, 30),
(170, 20), (160, 20)},
s0="prev"  @@{;"prev";**@{-}="prev"}
% arrows in X
\POS(18,18)\ar+(1,1), \POS(35,25)\ar+(2,0), \POS(35,45)\ar+(2,0), 
\POS(25,34)\ar+(0,2), \POS(45,34)\ar+(0,2), \POS(53,53)\ar+(1,1),
% arrows in Y
\POS(85,15)\ar+(1,1), \POS(90,32)\ar+(0,2), \POS(90,45)\ar+(0,2),
\POS(100,45)\ar+(0,2), \POS(110,50)\ar+(2,0), \POS(125,55)\ar+(1,1),
\POS(120,40)\ar+(0,2), \POS(120,25)\ar+(0,2), \POS(110,25)\ar+(0,2),
\POS(100,20)\ar+(2,0),
% arrows in Z
\POS(156,16)\ar+(1,1), \POS(175,35)\ar+(1,1), \POS(195,55)\ar+(1,1),
\POS(160,25)\ar+(0,2), \POS(170,25)\ar+(0,2), \POS(180,45)\ar+(0,2), \POS(190,45)\ar+(0,2), 
\endxy
\label{2.3.1} \end{array} \end{equation}

	A cell marked with a cross is not commutative, while the central cell in $ \uPi_1(\hat{Y}) $ commutes. 
In $ \uPi_1(\hat{X}) $ there are two arrows $ 0 \to 1$, in $ \uPi_1(\hat{Y}) $ there are three of them, in 
$ \uPi_1(\hat{Z}) $ four.

%%%
%%%
\section{Border flexible c-spaces and strict homotopies}\label{s3}

	We end by examining the relationship of border flexible c-spaces (defined in II.2.1(c)) with strict 
homotopies (see II.4.3(e)), expressed in Theorem \ref{3.1}. 

	As a consequence, the fundamental category of a border flexible c-space can be simply defined 
using c-paths up to homotopy with fixed endpoints (see \ref{3.2}). Its invariance up to strict homotopies is 
stated in Theorem \ref{3.3}.

	The importance of a simple construction, instead of the hybrid construction of Sections II.4 and II.5, is 
evident -- although it does not apply to essential c-spaces like the delayed intervals and the higher 
c-spheres, which are not border flexible (see II.2.2).

%%%
\subsection{Theorem {\rm (Border flexible c-spaces and homotopies)}}\label{3.1}
{\em
Let $Y$ be a border flexible c-space. Every strict homotopy $ \ph\c X \ti \cI \to Y $ is flexible.
}
\begin{proof}
We are given a c-map $ \ph\c X \times \cI \to Y $ which is constant on each fibre $ \{x\} \ti \cI $ at a flexible point of 
$ X$, and we have to prove that $ \ph $ is also a c-map $ X \ti \uI \to Y$.

	We take a c-path $ b = \lan a, h\ran \c \cI \to X \ti \uI$, where $ a\c x_0 \to x_1 $ is a c-path of $ X $ 
(between flexible points) and $ h\c t_0 \to t_1 $ is increasing in $ \uI$; we have to prove that $ \ph b $ is 
controlled in $ Y$.

	We insert a path $\, b_\al\c \cI \to X \ti \uI \,$ in each fibre of the cylinder at the endpoints 
$x_\al$ (for $\al = 0, 1$)
%
% FIGURE 3.1.1
    \begin{equation} 
\xy <.5mm, 0mm>:
%
% labels
(80,-10) *{X \ti \uI}; (10,-16) *{\sst{(x_0, 0)}}; (20,16) *{\sst{(x_1, 1)}}; 
(-5,-10) *{\sst{b_0}}; (12,5) *{\sst{b}}; (35,10) *{\sst{b_1}};
% bases of the cylinder  X \ti \uI
(-20,-10); (50,-10) **\crv{(-40,-30)&(70,-30)},
(-20,15); (50,15) **\crv{(-40,35)&(70,35)},
% curved path and its arrow
(0,-6); (30,6) **\crv{(5,5)&(25,-5)},
\POS(15,.1) \are+(4,.2), 
% dotted fibres, vertical paths, marks of endpoints
\POS(0,-15) \arlp+(0,30), \POS(30,-15) \arlp+(0,30), 
\POS(0,-15) \ar+(0,9), \POS(30,6) \ar+(0,9),
\POS(-1,-15) \arl+(2,0), \POS(29,15) \arl+(2,0),
\POS(-22.5,-13) \arl+(0,32), \POS(52.5,-13) \arl+(0,32), 
\endxy
    \label{3.1.1} \end{equation}
    \begin{equation*} \begin{array}{c}
b_0  =  \lan e_{x_0}, h_0 \ran \c \cI \to X \ti \uI,   \q   h_0\c 0 \to t_0,
\\[5pt]
b_1  =  \lan e_{x_1}, h_1\ran \c \cI \to X \ti \uI,   \q   h_1\c t_1 \to 1,
    \label{} \end{array} \end{equation*}
and we get a c-path $ b' = \lan a', h'\ran  = b_0 * b * b_1 $ in $ X \ti \uI $ which is controlled in $ X \ti \cI$, 
because $ h' $ is an increasing path $ 0 \to 1$.

	Now $ \ph b' $ is controlled in the border flexible c-space $ Y $ and each path $ \ph b_\al$ is constant 
(because $ \ph $ is a strict homotopy). It follows that the middle restriction $ \ph b $ is also controlled in $ Y$.
\end{proof}
%

%%%
\subsection{The border flexible case}\label{3.2}
As a particular case of the previous theorem, if the c-space $ X $ is border flexible, a general 2-path 
$\cI^2 \to X$ is always a hybrid 2-path $ \cI \ti \uI \to X$ (because $H$ is constant on the vertical edges of 
$\cI^2$).

	Therefore the restricted functor
    \begin{equation}
\uPi_1\c \cbfTop \to \Cat,
    \label{3.2.1} \end{equation}
can be equivalently defined using general 2-paths, based on the standard square $\cI^2$, instead of hybrid 
2-paths based on $\cI \ti \uI$.

	The restricted functor is still invariant up to flexible homotopies. But strict homotopies in $\cbfTop$ are 
always flexible, giving the following result.

%%%
\subsection{Theorem {\rm (Homotopy invariance, III)}}\label{3.3}
{\em
A strict homotopy $ \ph\c f \to g\c X \to Y $ of border flexible c-spaces induces the identity of the associated 
functors
    \begin{equation}
f_*  =  g_*\c \uPi_1(X) \to \uPi_1(Y).
    \label{3.3.1} \end{equation}
}%
\begin{proof}
By Theorem \ref{3.1}, $\ph$ is a strict flexible homotopy. Applying Theorem II.5.4(b), $\ph_*$ is the identity 
of $ f_* =  g_*$.	
\end{proof}
%

%%%
%%%

\vspace{5mm}
\noindent
    Marco Grandis\\
    Dipartimento di Matematica \\
    Universit\`a di Genova \\
    Via Dodecaneso 35 \\
    16146 - Genova, Italy \\
    grandis@dima.unige.it

\end{document}